\RequirePackage{fix-cm}
\documentclass[smallextended]{svjour3}       
\smartqed  
\usepackage{graphicx}
\usepackage{bm}
\usepackage{color}
\usepackage{ulem}

\usepackage{latexsym}
\def\qed{\hfill $\Box$}

\begin{document}

\title{Information Geometry Connecting Wasserstein Distance and Kullback-Leibler Divergence via the Entropy-Relaxed Transportation Problem 
}


\author{Shun-ichi Amari  \and Ryo Karakida \and Masafumi Oizumi }


\institute{Shun-ichi Amari \at
              2-1 Hirosawa, Wako-shi, Saitama, 351-0198, Japan \\
              Tel.: +81-48-467-9669 \\
              Fax: +81-48-467-9687 \\
              \email{amari@brain.riken.jp}           
              \and
           Ryo Karakida \at
           2-3-26 Aomi, Koto-ku, Tokyo, 135-0064, Japan \
          \and
          Masafumi Oizumi \at
          2-8-10 Toranomon, Minato-ku, Tokyo, 105-0001, Japan \\
}

\date{}

\maketitle

\begin{abstract} Two geometrical structures have been extensively studied for a manifold of probability distributions. One is based on the Fisher information metric, which is invariant under reversible transformations of random variables, while the other is based on the Wasserstein distance of optimal transportation, which reflects the structure of the distance between random variables. Here, we propose a new information-geometrical theory that is a unified framework connecting the Wasserstein distance and Kullback-Leibler (KL) divergence. We primarily considered a discrete case consisting of $n$ elements and studied the geometry of the probability simplex $S_{n-1}$, which is the set of all probability distributions over $n$ elements. The Wasserstein distance was introduced in $S_{n-1}$ by the optimal transportation of commodities from distribution ${\bm{p}}$ to distribution ${\bm{q}}$, where ${\bm{p}}$, ${\bm{q}} \in S_{n-1}$. We relaxed the optimal transportation by using entropy, which was introduced by Cuturi. The optimal solution was called the entropy-relaxed stochastic transportation plan. The entropy-relaxed optimal cost $C({\bm{p}}, {\bm{q}})$ was computationally much less demanding than the original Wasserstein distance but does not define a distance because it is not minimized at ${\bm{p}}={\bm{q}}$. To define a proper divergence while retaining the computational advantage, we first introduced a divergence function in the manifold $S_{n-1} \times S_{n-1}$ of optimal transportation plans. We fully explored the information geometry of the manifold of the optimal transportation plans and subsequently constructed a new one-parameter family of divergences in $S_{n-1}$ that are related to both the Wasserstein distance and the KL-divergence.

\end{abstract}

\section{Introduction}

Information geometry \cite{Amari2016} studies the properties of a manifold of probability distributions and is useful for various applications in statistics, machine learning, signal processing, and optimization. Two geometrical structures have been introduced from two distinct backgrounds. One is based on the invariance principle, where the geometry is invariant under reversible transformations of random variables. The Fisher information matrix, for example, is a unique invariant Riemannian metric from the invariance principle \cite{Amari2016,Chentsov1982,Rao1945}. Moreover, two dually coupled affine connections are used as invariant connections [1, 9], which are useful in various applications.

The other geometrical structure was introduced through the transportation problem, where one distribution of commodities is transported to another distribution. The minimum transportation cost defines a distance between the two distributions, which is called the Wasserstein, Kantorovich or earth-mover distance \cite{Santambrogio2015,Villani2013}. This structure provides a tool to study the geometry of distributions by taking the metric of the supporting manifold into account.

Let $X= \left\{1, \cdots, n \right\}$ be the support of a probability measure ${\bm{p}}$. The invariant geometry provides a structure that is invariant under permutations of elements of $X$ and results in an efficient estimator in statistical estimation. On the other hand, when we consider a picture over $n^2$ pixels $X= \left\{ (ij); i, j=1, \cdots, n \right\}$ and regard it as a distribution over $X$, the pixels have a proper distance structure in $X$. Spatially close pixels tend to take similar values. A permutation of $X$ destroys such a neighboring structure, suggesting that the invariance might not play a useful role. The Wasserstein distance takes such a structure into account and is therefore useful for problems with metric structure in support $X$ (see, e.g., \cite{Cuturi2013,Cuturi2014,Cuturi2016}).

An interesting question is how these two geometrical structures are related. While both are important in their own respects, it would be intriguing to construct a unified framework that connects the two. With this purpose in mind, we examined the discrete case over $n$ elements, such that a probability distribution is given by a probability vector ${\bm{p}}=(p,\cdots, p_n)$ in the probability simplex 
\begin{equation}
 S_{n-1} = \left\{ {\bm{p}}\; \left| \; p_i > 0, \;
 \sum p_i=1  \right.\right\}. 
\end{equation} We also consider Gaussian distributions over the one-dimensional real line $X$. 

Cuturi modified the transportation problem such that the cost is minimized under an entropy constraint \cite{Cuturi2013}. This is called the entropy-relaxed optimal translation problem and is computationally less demanding than the original transportation problem. In addition to the advantage in computational cost, Cuturi showed that the quasi-distance defined by the entropy-relaxed optimal solution yields superior results in many applications compared to the original Wasserstein distance and information-geometric divergences such as the KL divergence. 

We followed the entropy-relaxed framework that Cuturi et al. proposed \cite{Cuturi2013,Cuturi2014,Cuturi2016} and introduced a Lagrangian function, which is a linear combination of the transportation cost and entropy. Given a distribution ${\bm{p}}$ of commodity on the sender’s side and ${\bm{q}}$ on the receiver’s side, the constrained optimal transportation plan is the minimizer of the Lagrangian function. The minimum value $C({\bm{p}}, {\bm{q}})$ is a function of ${\bm{p}}$ and ${\bm{q}}$, which we called the Cuturi function. However, this does not define the distance between ${\bm{p}}$ and ${\bm{q}}$ because it is non-zero at ${\bm{p}} = {\bm{q}}$ and is not minimized when ${\bm{p}}$ = ${\bm{q}}$.

To define a proper distance-like function in $S_{n-1}$, we introduced a divergence between ${\bm{p}}$ and ${\bm{q}}$ derived from the optimal transportation plan. A divergence is a general metric concept that includes the square of a distance but is more flexible, allowing non-symmetricity between ${\bm{p}}$ and ${\bm{q}}$. A manifold equipped with a divergence yields a Riemannian metric with a pair of dual affine connections. Dually coupled geodesics are defined, which possess remarkable properties, generalizing the Riemannian geometry [1]. 

We studied the geometry of the entropy-relaxed optimal transportation plans within the framework of information geometry. They form an exponential family of probability distributions defined in the product manifold $S_{n-1} \times S_{n-1}$. Therefore, a dually flat structure was introduced. The $m$-flat coordinates are the expectation parameters $({\bm{p}}, {\bm{q}})$ and their dual, $e$-flat coordinates (canonical parameters) are $({\bm{\alpha}}, {\bm{\beta}})$, which are assigned from the minimax duality of nonlinear optimization problems. We can naturally defined a canonical divergence, that is the KL divergence $KL[({\bm{p}}, {\bm{q}}) : ({\bm{p}'}, {\bm{q}'})]$ between the two optimal transportation plans for $({\bm{p}}, {\bm{q}})$ and $({\bm{p}'}, {\bm{q}'})$, sending ${\bm{p}}$ to ${\bm{q}}$ and ${\bm{p}'}$ to ${\bm{q}'}$, respectively.

To define a divergence from ${\bm{p}}$ to ${\bm{q}}$ in $S_{n-1}$, we used the reference distribution $\bm{r}$. Given $\bm{r}$, we defined a divergence between ${\bm{p}}$ and ${\bm{q}}$ by $KL[({\bm{r}}, {\bm{p}}) : ({\bm{r}}, {\bm{q}})]$. There are a number of potential choices for $\bm{r}$: one is to use ${\bm{r}} = {\bm{p}}$ and another is to use the arithmetic or geometric mean of ${\bm{p}}$ and ${\bm{q}}$. These options yield one-parameter families of divergences connecting the Wasserstein distance and KL-divergence. Our work uncovers a novel direction for studying the geometry of a manifold of probability distributions by integrating the Wasserstein distance and KL divergence. 

\section{Entropy-Constrained Transportation Problem}

\begin{figure}[t]
 \begin{center}
  \includegraphics[width=60mm]{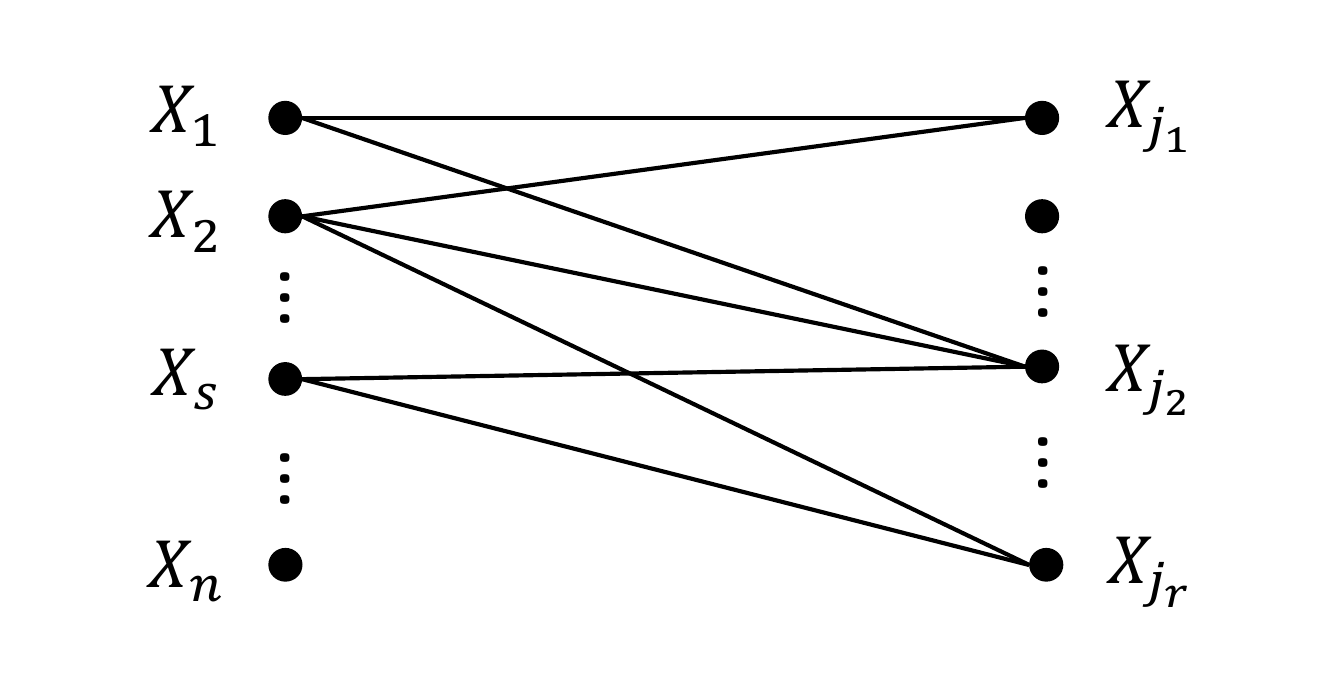}
 \end{center}
 \caption{Transportation from the sending terminals $X_S$ to the receiving terminals $X_R$}
 \label{fig1}
\end{figure}

Let us consider $n$ terminals $X= \left(X_1, \cdots, X_n \right)$, some of which, say $X_1, \cdots, X_s$, are sending terminals at which $p_1, \cdots, p_s$ ($p_i > 0$) amounts of commodities are stocked. At the other terminals, $X_{s+1}, \cdots, X_n$, no commodities are stocked ($p_i=0$). These are transported within $X$ such that $q_1, \cdots, q_r$ amounts are newly stored at the receiving terminals $X_{j_1}, \cdots, X_{j_r}$. There may be overlap in the sending and receiving terminals, $X_S = \left\{ X_1, \cdots, X_s \right\}$ and $X_R = \left\{X_{j_1}, \cdots, X_{j_r} \right\}$, including the case that $X_R=X_S=X$ (Fig. \ref{fig1}). We normalized the total amount of commodities to be equal to 1 so that ${\bm{p}}=\left(p_1, \cdots, p_s\right)$ and ${\bm{q}}= \left(q_1, \cdots, q_r \right)$ can be regarded as probability distributions in the probability simplex $S_{s-1}$ and $S_{r-1}$, respectively,
\begin{equation}
 \sum p_i=1, \quad \sum q_i=1, \quad p_i>0, \quad q_i>0.
\end{equation}
Let $S_{n-1}$ be the probability simplex over $X$.  Then $S_{s-1} \subset \bar{S}_{n-1}$, $S_{r-1} \subset \bar{S}_{n-1}$, where $\bar{S}_{n-1}$ is the closure of $S_{n-1}$,
\begin{equation}
 \bar{S}_{n-1} = \left\{ {\bm{r}}\; \left| \; r_i \ge 0, \;
 \sum r_i=1  \right.\right\}. 
\end{equation}

It should be noted that if some components of $\bm{p}$ and $\bm{q}$ are allowed to be 0, we do not need to treat $X_S$ and $X_R$ separately, i.e., we can consider both $X_S$ and $X_R$ to be equal to $X$. Under such a situation, we simply considered both $\bm{p}$ and $\bm{q}$ as elements of $\bar{S}_{n-1}$.

We considered a transportation plan ${\rm{\bf P}}=\left(P_{ij}\right)$ denoted by an $s \times r$ matrix, where $P_{ij} \ge 0$ is the amount of
commodity transported from $X_i \in X_S$ to $X_j \in X_R$. The plan ${\rm{\bf P}}$ was regarded as a (probability) distribution of commodities
flowing from $X_i$ to $X_j$, satisfying the sender and receiver’s conditions,
\begin{equation}
 \label{eq:am3}
 \sum_j P_{ij} = p_i, \quad \sum_i P_{ij}=q_j, \quad
 \sum_{ij} P_{ij} = 1.
\end{equation}
We denoted the set of ${\rm{\bf P}}$ satisfying Eq. (\ref{eq:am3}) as $U({\bm{p}}, {\bm{q}})$.

Let ${\rm{\bf M}}=\left(m_{ij}\right)$ be the cost matrix, where $m_{ij} \ge 0$ denotes the cost of transporting one unit of commodity from $X_i$
to $X_j$. We can interpret $m_{ij}$ as the distance between $X_i$ and $X_j$. The transportation cost of plan ${\rm{\bf P}}$ is 
\begin{equation}
  C({\rm{\bf P}}) = \langle {\rm{\bf M}}, {\rm{\bf P}} \rangle
 = \sum m_{ij} P_{ij}.
\end{equation}
The Wasserstein distance between ${\bm{p}}$ and ${\bm{q}}$ is the minimum cost of transporting commodities distributed by ${\bm{p}}$ at the sender’s to ${\bm{q}}$ at the receiver’s side,
\begin{equation}
  C_W({\bm{p}}, {\bm{q}}) = \mathop{\min}_{{\rm{\bf P}} 
  \subset U({\bm{p}}, {\bm{q}})}
 \langle {\rm{\bf M}}, {\rm{\bf P}} \rangle,
\end{equation}
where min is taken over all ${\rm{\bf P}}$ satisfying the constraints in Eq. (\ref{eq:am3}) \cite{Santambrogio2015,Villani2013}.

We considered the entropy of ${\rm{\bf P}}$, 
\begin{equation}
H({\rm{\bf P}})=-\sum P_{ij}
\log P_{ij}.  
\end{equation}
Given marginal distributions ${\bm{p}}$ and ${\bm{q}}$, the plan that maximizes the entropy is given by the direct product of ${\bm{p}}$ and ${\bm{q}}$,
\begin{equation}
 {\rm{\bf P}}_D = {\bm{p}} \otimes {\bm{q}} = \left(p_i q_j \right).
\end{equation}
This is because the entropy of ${\rm{\bf P}}_D$,
\begin{equation}
 H \left({\rm{\bf P}}_D \right) =
 -\sum {\rm{P}}_{Dij} \log {\rm{P}}_{Dij} =
 H({\bm{p}}) + H({\bm{q}}), 
\end{equation}
is the maximum among all possible ${\rm{\bf P}}$ belonging to $U({\bm{p}},
{\bm{q}})$, i.e.,
\begin{equation}
 H({\rm{\bf P}}) \le H({\bm{p}})+ H({\bm{q}}) = H( {\rm{\bf P}}_D ),
\end{equation}
where $H({\rm{\bf P}})$, $H({\bm{p}})$ and $H({\bm{q}})$ are the entropies of the respective distributions. 

We consider the constrained problem of searching for ${\rm{\bf P}}$ that minimizes $\langle {\rm{\bf M}}, {\rm{\bf P}} \rangle$ under the constraint $H({\rm{\bf P}}) \ge \mbox{const}$. This is equivalent to imposing the condition that ${\rm{\bf P}}$ lies within a KL-divergence ball centered at ${\rm{\bf P}}_D$,
\begin{equation}
 KL \left[{\rm{\bf P}}:{\rm{\bf P}}_D \right] \le d
\end{equation}
for constant $d$, because the KL-divergence from plan ${\rm{\bf P}}$ to ${\rm{\bf P}}_D$ is
\begin{equation}
 KL \left[{\rm{\bf P}}:{\rm{\bf P}}_D \right] =
 \sum P_{ij} \log \frac{P_{ij}}{p_i q_j} =
   -H({\rm{\bf P}}) + H({\bm{p}})+H({\bm{q}}).
\end{equation}
The entropy of ${\rm{\bf P}}$ increases within the ball as $d$ increases. Therefore, this is equivalent to the entropy constrained problem that minimizes a linear combination of the transportation cost $\langle {\rm{\bf M}}, {\rm{\bf P}} \rangle$ and entropy $H({\rm{\bf P}})$,
\begin{equation}
 \label{eq:am10}
 F_{\lambda}({\rm{\bf P}}) = 
 \langle {\rm{\bf M}}, {\rm{\bf P}} \rangle
 -\lambda H({\rm{\bf P}})
\end{equation}
for constant $\lambda$ \cite{Cuturi2013}. Here, $\lambda$ is a Lagrangian multiplier and $\lambda$ becomes smaller as $d$ becomes larger.

\section{Solution to the Entropy-Constrained Problem: Cuturi Function}

Let us fix $\lambda$ as a parameter controlling the magnitude of the entropy or the size of the KL-ball. When ${\rm{\bf P}}$ satisfies the constraints in Eq. (\ref{eq:am3}), minimization of Eq. (\ref{eq:am10}) is formulated in the Lagrangian form by using Lagrangian multipliers ${\alpha}_i$, ${\beta}_j$,
\begin{equation}
 \label{eq:am1220170111}
 L_{\lambda}({\rm{\bf P}}) = 
 \frac{1}{1+\lambda} \langle {\rm{\bf M}}, {\rm{\bf P}} \rangle
 -\frac{\lambda}{1+\lambda} H({\rm{\bf P}}) - \sum_{i, j} \left(\alpha_i + \beta_j
  \right)  P_{ij}.
\end{equation}
By differentiating Eq. (\ref{eq:am1220170111}) with respect to $P_{ij}$, we have
\begin{equation}
 \frac {1+\lambda}{\lambda}\frac{\partial}{\partial P_{ij}} L_{\lambda}
 ({\rm{\bf P}}) = \frac 1{\lambda} m_{ij}+ 
 \log P_{ij} -  \frac {1+\lambda}{\lambda} \left(\alpha_i+\beta_j \right) +1. 
\end{equation}

By setting the above derivatives equal to 0, we have the following solution,
\begin{equation}
 \label{eq:am1420170113}
  P_{ij} \propto \exp \left\{
 -\frac{m_{ij}}{\lambda}+  \frac {1+\lambda}{\lambda} \left(\alpha_i+\beta_j \right)  \right\}.
\end{equation}
Let us put
\begin{eqnarray}
 && K_{ij} = \exp \left\{-\frac{m_{ij}}{\lambda}\right\}, \label{eq:K} \\
 && a_i = \exp \left(\frac {1+\lambda}{\lambda} \alpha_i \right), \quad
 b_j = \exp \left( \frac {1+\lambda}{\lambda} \beta_j \right). 
\end{eqnarray}
Then, the optimal solution is written as
\begin{equation}
 P^{\ast}_{ij} = c a_i b_j K_{ij},
\end{equation}
where $a_i$ and $b_j$ are positive and correspond to the Lagrangian multipliers $\alpha_i$ and $\beta_j$ to be determined from the constraints (Eq. (\ref{eq:am3})). $c$ is the normalization constant. Since $r+s$ constraints (Eq. (\ref{eq:am3})) are not independent because of the conditions that $\sum p_i = 1$ and $\sum　q_j=1$, we can use $b_r=1$.  Further, we noted that $\mu {\bm{a}}$ and ${\bm{b}}/\mu$ yield the same answer for any $\mu>0$, where ${\bm{a}}=\left(a_i \right)$ and ${\bm{b}}= \left(b_j \right)$.  Therefore, the degrees of freedom of ${\bm{a}}$ and ${\bm{b}}$ are $s-1$ and $r-1$, respectively. We can choose ${\bm{a}}$ and ${\bm{b}}$ such that they satisfy
\begin{equation}
 \sum a_i = 1, \quad \sum b_j = 1.
\end{equation}
Then, $\bm{a}$ and $\bm{b}$ are included in $S_{s-1}$. We have the theorem below.

\theorem \upshape The optimal transportation plan ${\rm{\bf　P}}^{\ast}_{\lambda}$ is given by　
\begin{eqnarray}
 \label{eq:am19}
 P^{\ast}_{\lambda ij} &=& c a_i b_j K_{ij}, \\
 \label{eq:am20}
 c &=& \frac 1{\sum a_i b_j K_{ij}},
\end{eqnarray}
where two vectors ${\bm{a}}$ and ${\bm{b}}$ are determined from ${\bm{p}}$ and ${\bm{q}}$ using Eq. (\ref{eq:am3}).

We have a generalized cost function of transporting ${\bm{p}}$ to ${\bm{q}}$ based on the entropy-constrained optimal plan ${\rm{\bf P}}^{\ast}_{\lambda}({\bm{p}}, {\bm{q}})$: 
\begin{eqnarray}
 C_{\lambda}({\bm{p}}, {\bm{q}}) &=& \frac 1{1+\lambda} 
 \left<{\rm{\bf M}}, {\rm{\bf P}}^{\ast}_{\lambda} \right> 
 -\frac{\lambda}{1+\lambda} H \left({\rm{\bf
			       P}}^{\ast}_{\lambda}\right). \label{eq:Cuturi_cost}
\end{eqnarray}
We called it the Cuturi function because extensive studies have been conducted by Cuturi and colleagues \cite{Cuturi2013,Cuturi2014,Cuturi2016}. The function has been used in various applications as a measure of discrepancy between ${\bm{p}}$ and ${\bm{q}}$. The following theorem holds for the Cuturi function:

\theorem \upshape The Cuturi function $C_{\lambda} ({\bm{p}}, {\bm{q}})$ is a convex function of $({\bm{p}}, {\bm{q}})$.

\begin{proof} \upshape
Let ${\rm{\bf P}}^{\ast}_1$ and ${\rm{\bf P}}^{\ast}_2$ be the optimal solutions of transportation problems $\left({\bm{p}}_1, {\bm{q}}_1 \right)$ and $\left({\bm{p}}_2, {\bm{q}}_2 \right)$, respectively. For scalar $ 0  \le \nu \le 1$, we use
\begin{equation}
 \bar{\rm{\bf P}} = \nu {\rm{\bf P}}^{\ast}_1
  + (1-\nu){\rm{\bf P}}^{\ast}_2. 
\end{equation}
We have
\begin{eqnarray}
 \lefteqn{\nu C_{\lambda}\left({\bm{p}}_1:{\bm{q}}_1 \right) + 
 (1-\nu) C_{\lambda} \left({\bm{p}}_2 :{\bm{q}}_2 \right)} \nonumber \\
 && = \frac 1{(1+\lambda)} 
  \left\{ \nu \left< {\rm{\bf M}}, {\rm{\bf P}}^{\ast}_1 \right>
  +(1-\nu) \left<{\rm{\bf M}}, {\rm{\bf P}}^{\ast}_2 \right> \right\}
  -\frac{\lambda}{1+\lambda} 
  \left\{ \nu H \left({\rm{\bf P}}^{\ast}_1 \right) + (1-\nu)H
 \left({\rm{\bf P}}^{\ast}_2 \right)\right\} \nonumber \\
 && \ge \frac 1{(1+\lambda)} \left< {\rm{\bf M}}, 
  \bar{\rm{\bf P}} \right> - \frac{\lambda}{1+\lambda} H
 \left(\bar{\rm{\bf P}} \right),
\end{eqnarray}
because $H({\rm{\bf P}})$ is a concave function of ${\rm{\bf P}}$. We further have
\begin{eqnarray}
 && \frac 1{(1+\lambda)} \langle {\rm{\bf M}},\bar{\rm{\bf P}} \rangle
   - \frac{\lambda}{1+\lambda} H
 \left(\bar{\rm{\bf P}}\right) \ge {\mathop{\min}_{\rm{\bf P}}}
 \left\{ \frac 1{(1+\lambda)} \langle {\rm{\bf M}}, {\rm{\bf P}}
  \rangle
  - \frac{\lambda}{1+\lambda} H({\rm{\bf P}}) \right\} \nonumber \\
 &&\qquad = C_{\lambda} \left\{\nu {\bf p}_1 +
			       (1-\nu){\bf p}_2, \nu {\bm{q}}_1 + (1-\nu)
 {\bm{q}}_2 \right\},
\end{eqnarray}
since the minimum is taken for ${\rm{\bf P}}$ transporting commodities from $\nu{\bm{p}}_1 + (1-\nu){\bm{p}}_2$ to $\nu {\bm{q}}_1+
(1-\nu){\bm{q}}_2$. Hence, the convexity of $C_{\lambda}$ is proven. \qed
\end{proof}

When $\lambda \rightarrow 0$, it converges to the original Wasserstein distance $C_W({\bm{p}}, {\bm{q}})$. However, it does not satisfy important requirements for ``distance''. When ${\bm{p}}={\bm{q}}$, $C_{\lambda}$ is not equal to
0 and does not take the minimum value, i.e., there are some ${\bm{q}}$ ($\ne {\bm{p}}$) that yield smaller $C_{\lambda}$ than ${\bm{q}} =\bm{p}$:
\begin{eqnarray}
 C_{\lambda}({\bm{p}}, {\bm{p}}) &>& C_{\lambda}({\bm{p}}, {\bm{q}})
\end{eqnarray}

\section{Geometry of Optimal Transportation Plans}
We first showed that a set of optimal transportation plans forms an exponential family embedded within the manifold of all transportation plans. Then, we studied the invariant geometry induced within these plans. A transportation plan ${\rm{\bf P}}$ is a probability distribution over branches $(i, j)$ connecting terminals of $X_i \in X_S$ and $X_j \in X_R$. Let $x$ denote branches $(i, j)$. We used the delta function $\delta_{ij}(x)$, which is 1 when $x$ is $(i, j)$ and 0 otherwise. Then, ${\rm{\bf P}}$ is written as a probability distribution of the random variable $x$,
\begin{equation}
 P(x) = \sum_{i,j} P_{ij} \delta_{ij}(x).
\end{equation}
By introducing new parameters
\begin{equation}
 \theta^{ij} = \log \frac{P_{ij}}{P_{sr}}, \quad
 {\bm{\theta}} = \left(\theta^{ij}\right),
\end{equation}
it is rewritten in parameterized form as
\begin{equation}
 P(x, {\bm{\theta}}) = \exp \left\{
 \sum_{i, j} \theta^{ij} \delta_{ij}(x) + \log P_{sr}
 \right\}.
\end{equation}
This shows that the set of transportation plans is an exponential family, where $\theta^{ij}$ are the canonical parameters and $\eta_{ij}=P_{ij}$ are the expectation parameters. They form an $\left(sr-1\right)$-dimensional manifold denoted by $S_{TP}$, because $\theta^{sr}=0$.

The transportation problem is related to various problems in information theory such as the rate-distortion theory. We provide detailed studies on the transportation plans in the information-geometric framework in Section 7, but here we introduce the manifold of the optimal transportation plans, which are determined by the sender’s and receiver’s probability distributions $\bm{p}$ and $\bm{q}$.

The optimal transportation plan specified by $({\bm{\alpha}}, {\bm{\beta}})$ in Eq. (\ref{eq:am1420170113}) is written as
\begin{equation}
 P(x, {\bm{\alpha}}, {\bm{\beta}}) = \exp
 \left[ \sum_{i, j} \left\{
 \frac {1+\lambda}{\lambda} \left(\alpha_i+\beta_j \right)
  -\frac{m_{ij}}{\lambda}\right\}
  \delta_{ij}(x) -  \frac {1+\lambda}{\lambda} \psi  \right].
\end{equation}
 The notation $\psi$ is a normalization factor called the potential function which is defined by
\begin{equation}
 \label{eq:am2920170124}
 \psi({\bm{\alpha}}, {\bm{\beta}}) =  -\frac{\lambda}{1+\lambda}\log c,
\end{equation}
 where $c$ is calculated by taking the summation over all of $x$,  
 \begin{equation}
 c = \sum_{x \in (X_S,X_R)} \exp
 \left[ \sum_{i, j} \left\{
 \frac {1+\lambda}{\lambda} \left(\alpha_i+\beta_j \right)
  -\frac{m_{ij}}{\lambda}\right\}
  \delta_{ij}(x) \right].
\end{equation}

This corresponds to the free energy in physics. 
 By using
\begin{equation}
 \label{eq:am3020170124}
 \theta^{ij} =  \frac {1+\lambda}{\lambda}  \left(\alpha_i+\beta_j \right)
  -\frac{m_{ij}}{\lambda},
\end{equation}
we see that the set $S_{OTP}$ of the optimal transformation plans is a submanifold of $S_{TP}$. Because Eq. (\ref{eq:am3020170124}) is linear in ${\bm{\alpha}}$ and ${\bm{\beta}}$, $S_{OTP}$ itself is an exponential family, where the canonical parameters are $({\bm{\alpha}}, {\bm{\beta}})$ and the expectation parameters are $({\bm{p}}, {\bm{q}}) \in S_{s-1} \times S_{r-1}$. This is confirmed by
\begin{eqnarray}
 {\rm{E}} \left[ \sum_j \delta_{ij} (x)\right] &=& p_i, \\
 {\rm{E}} \left[ \sum_i \delta_{ij}(x)\right] &=& q_j,
\end{eqnarray}
where ${\rm{E}}$ denotes the expectation. Because of ${\bm{p}} \in S_{s-1}$ and ${\bm{q}} \in S_{r-1}$, $S_{OPT}$ is a $(r+s-2)$-dimensional dually flat manifold, We can use $\alpha_s=\beta_r=0$ without loss of generality, which corresponds to using $a_s=b_r=1$ instead of the normalization $\sum a_i = \sum b_j = 1$ of ${\bm{a}}$ and ${\bm{b}}$.

In a dually flat manifold, the dual potential function $\varphi_{\lambda}$ is given from the potential function $\psi_{\lambda}$ as its Legendre dual, which is given by

\begin{eqnarray}
 \varphi_{\lambda} ({\bm{p}}, {\bm{q}}) &=&  {\bm{p}} \cdot {\bm{\alpha}} + {\bm{q}} \cdot
  {\bm{\beta}} - \psi_{\lambda}({\bm{\alpha}}, {\bm{\beta}}). \label{eq:rescaleC}
\end{eqnarray}
When we use new notations ${\bm{\eta}}= ({\bm{p}}, {\bm{q}})^T$,
${\bm{\theta}}= ({\bm{\alpha}}, {\bm{\beta}})^T$, we have 
\begin{equation}
 \psi_{\lambda}({\bm{\theta}}) + \varphi_{\lambda}({\bm{\eta}})
 = {\bm{\theta}} \cdot {\bm{\eta}},
\end{equation}
which is the Legendre relationship between ${\bm{\theta}}$ and ${\bm{\eta}}$, we have the following theorem:

\theorem \upshape The dual potential $\varphi_{\lambda}$ is equivalent to the Cuturi function $C_{\lambda}$.

\begin{proof} \upshape
Direct calculation of Eq. (\ref{eq:rescaleC}) gives 
\begin{eqnarray}
\varphi_{\lambda} ({\bm{p}}, {\bm{q}}) &=& 
{\bm{p}} \cdot {\bm{\alpha}} + {\bm{q}} \cdot
  {\bm{\beta}} - \psi_{\lambda}({\bm{\alpha}}, {\bm{\beta}}) \nonumber \\
&=& \frac{1}{1+\lambda}
 \left<{\rm{\bf M}}, {\rm{\bf P}} \right> +
 \sum_{i, j} P_{ij}
 \left\{ \left(\alpha_i + \beta_j \right)
  - \frac{1}{1+\lambda} m_{ij} - \psi_{\lambda}\right\} \nonumber \\
  &=&  \frac{1}{1+\lambda} \left<{\rm{\bf M}}, {\rm{\bf P}} \right> +
  \frac{\lambda}{1+\lambda} \sum_{i, j} P_{ij}
 \left ( \log  a_i + \log b_j - \frac{m_{ij}}{\lambda} + \log c \right ) \nonumber \\
  &=& C_{\lambda}({\bm{p}}, {\bm{q}}). 
\end{eqnarray}
\qed  \end{proof}

We summarize the Legendre relationship below.
\theorem \upshape The dual potential function $\varphi_{\lambda}$ (Cuturi function) and potential function (free energy, cumulant generating function) $\psi_{\lambda}$ of the exponential family $S_{OPT}$ are both convex, connected by the Legendre transformation,
\begin{equation}
   {\bm{\theta}} = \nabla_{\bm{\eta}} \varphi_{\lambda}({\bm{\eta}}), 
  \quad {\bm{\eta}} = \nabla_{\bm{\theta}}
  \psi_{\lambda}({\bm{\theta}}), 
\end{equation}
or
\begin{eqnarray}
  & {\bm{\alpha}} = \nabla_{\bm{p}} \varphi_{\lambda} ({\bm{p}},
   {\bm{q}}), &
 {\bm{\beta}} = \nabla_{\bm{q}} \varphi_{\lambda} ({\bm{p}}, {\bm{q}}),
 \\
 & {\bm{p}} = \nabla_{\bm{\alpha}} \psi_{\lambda}({\bm{\alpha}},
  {\bm{\beta}}), & {\bm{q}} = \nabla_{\bm{\beta}} \psi_{\lambda}
  ({\bm{\alpha}}, {\bm{\beta}}).
\end{eqnarray}

Since $S_{OPT}$ is dually flat, we can introduce a Riemannian metric and cubic tensor. 
The Riemannian metric ${\rm{\bf G}}_{\lambda}$ is given to $S_{s-1} \times S_{r-1}$ by
\begin{equation}
 {\rm{\bf G}}_{\lambda} = \nabla_{\bm{\eta}} \nabla_{\bm{\eta}}
 \varphi_{\lambda}({\bm{\eta}})
\end{equation}
in the ${\bm{\eta}}$-coordinate system $({\bm{p}}, {\bm{q}})$. Its inverse is
\begin{equation}
 \label{eq:am38-20170308}
 {\rm{\bf G}}^{-1}_{\lambda} = \nabla_{\bm{\theta}}
  \nabla_{\bm{\theta}} \psi_{\lambda}({\bm{\theta}}).
\end{equation}
Calculating Eq. (\ref{eq:am38-20170308}) carefully, we have the following theorem:

\theorem \upshape
The Fisher information matrix ${\rm{\bf G}}^{-1}_{\lambda}$ in the
${\bm{\theta}}$-coordinate system is given by
\begin{equation}
 {\rm{\bf G}}^{-1}_{\lambda} = 
 \left[
 \begin{array}{c|c}
   p_i \delta_{ij}-p_i p_j & P_{ij}-p_i q_j  \\
   \hline
    P_{ij}-p_i q_j & q_i \delta_{ij}-q_i q_j
  \end{array}
 \right].
\end{equation}

\begin{description}
\item[Remark 1.]  The ${\bm{p}}$-part and ${\bm{q}}$-part of ${\rm{\bf G}}^{-1}_{\lambda}$ are equal to the corresponding Fisher information in $S_{s-1}$ and $S_{r-1}$ in the $e$-coordinate systems. 
\item[Remark 2.]  The ${\bm{p}}$-part and the ${\bm{q}}$-part of ${\rm{\bf G}}_{\lambda}$ are not equal to the corresponding Fisher information in the $m$-coordinate system. This is because $({\bm{p}}, {\bm{q}})$-part of ${\rm{\bf G}}$ is not 0.
\end{description}
We can similarly calculate the cubic tensor,
\begin{equation}
 {\rm{\bf T}} = \nabla \nabla \nabla \psi_{\lambda}
\end{equation}
but we have not shown the results here. 

From the Legendre pair of convex functions $\varphi_{\lambda}$ and $\psi_{\lambda}$, we can also introduce the canonical divergence between two transportation problems $({\bm{p}}, {\bm{q}})$ and $({\bm{p}}', {\bm{q}}')$, 
\begin{equation}
D_{\lambda} \left[({\bm{p}}, {\bm{q}}) : \left({\bm{p}}', {\bm{q}}' \right)\right]  = \psi_{\lambda} ({\bm{\alpha}}, {\bm{\beta}})  + \varphi_{\lambda} ({\bm{p}}', {\bm{q}}')  - {\bm{\alpha}} \cdot {\bm{p}}'  - {\bm{\beta}} \cdot {\bm{q}}' 
\end{equation}
where $({\bm{\alpha}}, {\bm{\beta}})$ corresponds to $({\bm{p}}, {\bm{q}})$. \
This is the KL-divergence between the two optimal transportation plans,
\begin{equation}
D_{\lambda} \left[({\bm{p}}, {\bm{q}}) : \left({\bm{p}}', {\bm{q}}' \right)\right] 
= KL [P_{\lambda}({\bm{p}},{\bm{q}}): P_{\lambda}({\bm{p'}},{\bm{q'}})]. \label{eq:canonical}
\end{equation}

\section{$\lambda$-Divergences in $S_{n-1}$}

\subsection{Derivation of $\lambda$-divergences}
We defined the divergence between $\bm{p} \in S_{n-1}$ and $\bm{q} \in S_{n-1}$ using the canonical divergence in the set $S_{OTP}$ of the optimal transportation plans (Eq. (\ref{eq:canonical})). For the sake of simplicity, we hereafter only studied the case $X_S = X_R = X$.　 We introduce a reference distribution $\bm{r} \in S_{n-1} $ and defined the $\bm{r}$-referenced divergence between $\bm{p}$ and $\bm{q}$ by 

\begin{equation}
 D_{{\bm{r}}, \lambda}[{\bm{p}}:{\bm{q}}] = \gamma_{\lambda} KL 
 \left[{\rm{\bf P}}_{\lambda}^*({\bm{r}}, {\bm{p}}):
 {\rm{\bf P}}^*_{\lambda}({\bm{r}}, {\bm{q}})\right],
\end{equation}
where $\gamma_{\lambda}$ is a scaling factor, which we discuss later, and ${\rm{\bf P}}^{\ast}_{\lambda}({\bm{r}}, {\bm{p}})$ is the optimal transportation plan from ${\bm{r}}$ to ${\bm{p}}$.

There are various ways of choosing a reference distribution $\bm{r}$. We first considered the simple choice of $\bm{r}=\bm{p}$, yielding the following $\lambda$-divergence:
\begin{equation}
 D_{\lambda}[{\bm{p}}:{\bm{q}}] = \gamma_{\lambda} KL 
 \left[{\rm{\bf P}}^*_{\lambda}({\bm{p}}, {\bm{p}}):
 {\rm{\bf P}}^*_{\lambda}({\bm{p}}, {\bm{q}})\right].
\end{equation}

\theorem \upshape $D_{\lambda}[{\bm{p}}: {\bm{q}}]$ with the scaling factor $\gamma_{\lambda}= \frac{\lambda}{1+\lambda}$ is given by
\begin{equation}
 D_{\lambda}[{\bm{p}}:{\bm{q}}] =  C_{\lambda}({\bm{p}}, {\bm{p}})
  -C_{\lambda}({\bm{p}}, {\bm{q}}) -\nabla_{\bm{q}}C_{\lambda}({\bm{p}}, {\bm{q}}) \cdot  ({\bm{p}}-{\bm{q}}), \label{eq:Bregman_C}
\end{equation}
which is constructed from the Cuturi function.

\begin{proof}
The optimal transportation plans are rewritten by the $\bf{\theta}$ coordinates in the form 
\begin{eqnarray}
&& \frac{\lambda}{1+\lambda} \log {\rm{\bf P}}^*_{\lambda}({\bm{p}}, {\bm{p}})_{ij}
= \alpha'_i + \beta'_j  -\frac{m_{ij}}{\lambda}-\psi'_{\lambda}, \\
&& \frac{\lambda}{1+\lambda} \log {\rm{\bf P}}^*_{\lambda}({\bm{p}}, {\bm{q}})_{ij}
= \alpha_i + \beta_j  -\frac{m_{ij}}{\lambda}-\psi_{\lambda}.
\end{eqnarray}
Then, we have
\begin{eqnarray}
 && D_{\lambda}[{\bm{p}}:{\bm{q}}] =   {\bm{p}} \cdot {\bm{\alpha}}' + {\bm{p}} \cdot
  {\bm{\beta}}' - \psi_{\lambda}' - {\bm{p}} \cdot {\bm{\alpha}} - {\bm{q}} \cdot
  {\bm{\beta}} - \psi_{\lambda} - ( {\bm{p}}- {\bm{q}}) \cdot {\bm{\beta}}  \nonumber \\
 && \mbox{\qquad\qquad}= \varphi_{\lambda}({\bm{p}}, {\bm{p}})
  -\varphi_{\lambda}({\bm{p}}, {\bm{q}}) -
  \nabla_{\bm{q}} \varphi_{\lambda}({\bm{p}}, {\bm{q}}) \cdot
  ({\bm{p}}-{\bm{q}}).
\end{eqnarray}
Since we showed that $\varphi_{\lambda}=C_{\lambda}$ in Theorem 3, we obtain Eq. (\ref{eq:Bregman_C}). 
\qed
\end{proof}

This is a divergence function satisfying $D_{\lambda}[{\bm{p}}:{\bm{q}}] \ge 0$, with equality when and only when ${\bm{p}}={\bm{q}}$. However, it is not a canonical divergence of a dually flat manifold. The Bregman divergence derived from a convex function $\tilde{\varphi}({\bm{p}})$ is given by
\begin{equation}
\tilde{D}_{\lambda} [{\bm{p}}:{\bm{q}}] = \tilde{\varphi}({\bm{p}}) - \tilde{\varphi}({\bm{q}}) - \nabla_{{\bm{p}}} \tilde{\varphi}({\bm{q}}) \cdot ({\bm{p}}- {\bm{q}}).
\end{equation}
This is different from Eq. (\ref{eq:Bregman_C}), which is derived from $\varphi_{\lambda}({\bm{p}},{\bm{q}})$. Thus, we call $D_{\lambda}[{\bm{p}}:{\bm{q}}]$ Bregman-like divergence.

In the extremes of $\lambda$, the proposed divergence $D_{\lambda}[{\bm{p}}:{\bm{q}}]$ is related to the KL-divergence and Wasserstein distance in the following sense:
\begin{enumerate}
 \item When $\lambda \rightarrow \infty$, $D_{\lambda}$ converges to $KL[{\bm{p}}:{\bm{q}}]$. This is because ${\rm{\bf P}}^*$ converges to ${\bm{p}} \otimes {\bm{q}}$ in the limit and we easily have
\begin{equation}
         KL[ {\bm{p}} \otimes {\bm{p}}: {\bm{p}} \otimes {\bm{q}}] = KL[\bm{p} : \bm{q}].
\end{equation}

\item When $\lambda \rightarrow 0$, $D_{\lambda}$ converges to 0, because $KL 
 \left[{\rm{\bf P}}_{0}^*({\bm{p}}, {\bm{p}}):
 {\rm{\bf P}}^*_{0}({\bm{p}}, {\bm{q}})\right]$ takes a finite value (see Example 1 in the next section). 
$C_\lambda=\varphi_\lambda$ is not differentiable when $\lambda =0$. Hence,  we cannot construct
the Bregman-like divergence from $C_0$ (Eq. (\ref{eq:Bregman_C})). This suggests that it is preferable to use a scaling factor other than $\gamma_{\lambda}= \lambda/(1+\lambda)$ when $\lambda$ is small. 
 \end{enumerate}

Since we have,
\begin{equation}
 \partial_{q_j} D_{\lambda }[{\bm{p}}:{\bm{q}}] =
 -\sum_i p_i \partial_{q_j} \log c a_i b_i =
 -\sum_i p_i \partial_{q_j} \log P_{ii},
\end{equation}
the Fisher information derived from $D_{\lambda}[{\bm{p}}:{\bm{q}}]$ is
\begin{equation}
 {\rm{\bf G}}_{\lambda} = \left. -\partial^2_{q_i q_j}
 D_{\lambda}[{\bm{p}}:{\bm{q}}]　\right|_{\bm{q}=\bm{p}} = \sum_k
 p_k \partial^2_{q_i q_j} \log P_{kk}.
\end{equation}

\subsection{Other choices of reference distribution $\bm{r}$} 
We can consider other choices of the reference distribution $\bm{r}$. One option is choosing $\bm{r}$, which minimizes the KL-divergence.
\begin{equation}
 \tilde{D}_{\lambda}[{\bm{p}}:{\bm{q}}] =  \gamma_{\lambda} 
 {\mathop{\min}_{\bm{r}}} KL
 \left[{\rm{\bf P}}_{\lambda}({\bm{p}}, {\bm{r}}):
 {\rm{\bf P}}_{\lambda}({\bm{q}}, {\bm{r}})\right].
\end{equation}
However, obtaining the minimizer $\bm{r}$ is not computationally easy. Thus, we can simply replace the optimal $\bm{r}$ with the arithmetic mean or geometric mean of $\bm{p}$ and $\bm{q}$. The arithmetic mean is given by the $m$-mixture midpoint of ${\bm{p}}$ and ${\bm{q}}$, 
\begin{equation}
 {\bm{r}} = \frac 12 ({\bm{p}}+{\bm{q}}).
\end{equation}
The geometric mean is given by the $e$-midpoint of ${\bm{p}}$ and ${\bm{q}}$,
\begin{equation}
 {\bm{r}} = c( \sqrt{p_i q_i}).
\end{equation}

\subsection{Examples of $\lambda$-Divergence}
Below, we consider the case where $\bm{r}=\bm{p}$. We show two simple examples, where $D_{\lambda}({\bm{p}}, {\bm{q}})$ can be analytically computed.

\ \

\noindent \textbf{Example 1}

Let $n=2$ and
\begin{equation}
 m_{ii} = 0, \quad m_{ij} = 1 \quad (i \ne j).
\end{equation}
We use $a_2=b_2=1$ for normalization,
\begin{eqnarray}
 P_{ij} &=& c a_i b_j K_{ij}, \\
 K_{ij} &=& \exp \left\{ -\frac{m_{ij}}{\lambda}\right\}
 = \left[
  \begin{array}{cc}
   1 & \varepsilon \\
   \varepsilon & 1
  \end{array}
 \right],\\
 \varepsilon &=& \exp \left\{ -\frac 1{\lambda}\right\}.
\end{eqnarray}
Note that $\varepsilon \rightarrow 0$ as $\lambda \rightarrow 0$.

When $\lambda>0$,
the receiver conditions require 
\begin{eqnarray}
 && cab + ca \varepsilon = p, \\
 && cab + cb \varepsilon = q, 
\end{eqnarray}
where we use $a=a_1$, $b=b_1$ and
\begin{equation}
 c = \frac 1{ab+ \varepsilon(a+b)+1}.
\end{equation}
Solving the above equations, we have
\begin{eqnarray}
 a &=& \frac{z-(q\!-\!p)/\varepsilon}{2(1\!-p)},  \\
 b &=& \frac{z+(q\!-\!p)/\varepsilon}{2(1\!-q)}, 
\end{eqnarray}
where 
\begin{eqnarray*}
 z &=&  - \varepsilon(1\!-p\!-q)+
  \sqrt{(q\!-\!p)^2/\varepsilon^2+\varepsilon^2(1\!-p\!-q)^2+2p(1\!-\!p)
  +2q(1\!-q)}.
\end{eqnarray*}

We can show $D_{\lambda}[{\bm{p}}:{\bm{q}}]$ explicitly by using the solution,
although it is complicated.

When $\lambda=0$, we easily have
\begin{equation}
 C_0(p, q) = |p-q|,
\end{equation}
where ${\bm{p}}=(p, 1-p)$ and ${\bm{q}}=(q, 1-q)$.  $C_0(p, q)$ is piecewise linear, and cannot be used to construct a Bregman-like divergence. However, we can calculate the limiting case of $\lambda \to 0$ because the optimal transportation plans ${\rm{\bf P}}^*$ where $\lambda$ is small are directly calculated by minimizing $C_\lambda (\bm{p},\bm{q})$ as 
\begin{eqnarray}
{\rm{\bf P}}^*_{\lambda}(\bm{p}, \bm{p}) &=& \left[
 \begin{array}{cc}
   p & 0  \\
    0 &1-p
  \end{array}
 \right] + \left[
 \begin{array}{cc}
   - \varepsilon &  \varepsilon   \\
     \varepsilon  & -  \varepsilon 
  \end{array}
 \right], \\
 {\rm{\bf P}}^*_{\lambda}(\bm{p}, \bm{q}) &=& \left[
 \begin{array}{cc}
   p & 0  \\
    q-p &1-q
  \end{array}
 \right] + \left[
 \begin{array}{cc}
   - \varepsilon^2 &   \varepsilon^2   \\
     \varepsilon^2  & -  \varepsilon^2 
  \end{array}
 \right].
\end{eqnarray} 
where we set $q>p$.
The limit of $KL$ divergence is given by 
\begin{eqnarray}
  \lim_{\lambda \to 0} KL[{\rm{\bf P}}^*_{\lambda}(\bm{p}, \bm{p}) :{\rm{\bf P}}^*_{\lambda}(\bm{p}, \bm{q}) ] = \left \{ \begin{array}{ll}
    p \log \frac{p}{q} & (p \geq q),   \\
   (1 - p) \log \frac{1 - p}{1 - q}  &  (p<q). 
  \end{array} \right.
\end{eqnarray}

In the general case of $n \geq 2$, the optimal transportation plan is 
${\rm{\bf P}}^*_{0}(\bm{p}, \bm{p}) = (p_i \delta_{ij})$. 
The diagonal parts of the optimal ${\rm{\bf P}}^*_{0}(\bm{p}, \bm{q}) $ are $\min \{p_i,q_i\}$ when $m_{ii}=0, \ \ m_{ij}>0 \ \ (i \neq j)$. Thus, the $KL $ divergence is given by 
\begin{equation}
KL[{\rm{\bf P}}^*_{0}(\bm{p}, \bm{p}) :{\rm{\bf P}}^*_{0}(\bm{p}, \bm{q}) ] = \sum_{i; p_i>q_i} p_i \log \frac{p_i}{q_i}.
\end{equation}
Remark that when $\lambda \rightarrow \infty$, 
\begin{equation}
  \lim_{\lambda \to \infty} KL[{\rm{\bf P}}^*_{\lambda}(\bm{p}, \bm{p}) :{\rm{\bf P}}^*_{\lambda}(\bm{p}, \bm{q}) ] = \sum_{i} p_i \log \frac{p_i}{q_i}.
\end{equation}

\ \ 

\noindent
{\textbf{Example 2}}

We take a family of Gaussian distributions $N \left(\mu, \sigma^2 \right)$, 
\begin{equation}
 p \left(x \;;\; \mu, \sigma^2 \right) =
 \frac 1{\sqrt{2 \pi} \sigma} \exp 
 \left\{ -\frac{(x-\mu)^2}{2 \sigma^2}\right\}
\end{equation}
on the real line $X= \left\{x \right\}$, extending the discrete case to the continuous case. We transport $p \left(x \;;\; \mu_p, \sigma^2_p \right)$ to $q \left(x\;;\; \mu_q, {\sigma}^2_q \right)$, where the transportation cost is
\begin{equation}
 m(x, y) = |x-y|^2.
\end{equation}
Then, we have
\begin{equation}
 K(x, y) = \exp \left\{ -\frac{(x-y)^2}{2 \lambda^2}\right\},
\end{equation}
where we use $2 \lambda^2$ instead of previous $\lambda$ for the sake of convenience.

The optimal transportation plan is written as
\begin{equation}
 P^{\ast}(x, y) = ca(x)b(y)K(x, y), \label{GaussianP}
\end{equation}
where $a$ and $b$ are determined from
\begin{eqnarray}
 \int ca(x)b(y)K(x, y)dy &=& p(x), \\
 \int ca(x)b(y)K(x, y)dx &=& q(x).
\end{eqnarray}
The solutions are given in the Gaussian framework, $x \sim N
\left(\tilde{\mu}, \tilde{\sigma}^2 \right)$, $y \sim N
\left(\tilde{\mu}', \tilde{\sigma}'^2 \right)$.
As derived in Appendix A,
the optimal cost and divergence are as follows:
\begin{eqnarray}
&&  C_{\lambda}(p,q) 
   = \frac 1{1+\lambda}
 \Biggl[ \left(\mu_p-\mu_q \right)^2 + \sigma^2_p +
 \sigma^2_q + \frac{\lambda}2 (1-\sqrt{1+X}) \nonumber \\
  && \mbox{\qquad\qquad\qquad}-\lambda
 \left \{ \log \sigma_p \sigma_q +\frac 12 \log 8 \pi^2e^2 - \frac 12 \log
 \left(1+\sqrt{1+X}\right)\right \} \Biggr], \label{GaussianC} \\
 && D_{\lambda} \left[p : q \right] = \gamma_{\lambda}
 \Biggl[ \frac 12 \left(\sqrt{1+X}-\sqrt{1+X_p}\right)
 + \log \frac{\sigma_q}{\sigma_p} + \frac 12 \log
 \frac{1+\sqrt{1+X_p}}{1+\sqrt{1+X}} \nonumber \\
 &&\mbox{\qquad\qquad\qquad\qquad} + \frac{1+\sqrt{1+X}}4
 \left\{ \frac{\left(\mu_p-\mu_q \right)^2}{\sigma^2_q}+
  \frac{\sigma^2_p}{\sigma^2_q} - 1 \right\}\Biggr],  \label{GaussianD} \\
  &&\mbox{\qquad\qquad}
   \mbox{where} \quad X= \frac{16 \sigma^2_p \sigma^2_q}{\lambda^2} \quad
 X_p = \frac{16 \sigma^4_p}{\lambda^2}. \nonumber
\end{eqnarray}

Note that $D_{\lambda}=KL \left[{\rm{\bf P}}^*_{\lambda}({\bm{p}}, {\bm{p}}): {\rm{\bf P}}^*_{\lambda}({\bm{p}}, {\bm{q}})\right]$ diverges to infinity in the limit of $\lambda \to 0$ because the support of the optimal transport ${\bf {P}}^*_{\lambda}({\bm{p}}, {\bm{q}})$ reduces to a 1-dimensional subspace. To prevent $D_{\lambda}$ from diverging and to make it finite, we set the scaling factor as $\gamma_{\lambda} = \frac{\lambda}{1+\lambda}$. In this case, $D_{\lambda}$ is equivalent to the Bregman-like divergence of the Cuturi function as shown in Theorem 6. With this scaling factor $\gamma_{\lambda}$, $D_{\lambda}$ in the limits of $\lambda \to \infty$ and $\lambda \to 0$ is given by
\begin{equation}
\lim_{\lambda \to \infty } D_{\lambda} = \frac{1}{2} \left \{ \frac{(\mu_p-\mu_q)^2}{\sigma_q^2} + \frac{\sigma_p^2}{\sigma_q^2} - 1 \right \} + \log \frac{\sigma_q}{\sigma_p} = KL[p : q],
\end{equation}
\begin{equation}
\lim_{\lambda \to 0 } D_{\lambda} = \frac{\sigma_p}{\sigma_q} (\mu_p-\mu_q)^2 + \frac{\sigma_p}{\sigma_q} (\sigma_p - \sigma_q)^2.
\end{equation}

\section{Applications of $\lambda $-Divergence}

\subsection{Cluster center (barycenter)}
Let ${\bm{q}}_1, \cdots, {\bm{q}}_k$ be $k$ distributions in $S_{n-1}$.
Its $\lambda$-center is defined by ${\bm{p}}^{\ast}$, which minimizes the
average of $\lambda$-divergences from ${\bm{q}}_i$ to ${\bm{p}} \in
S_{n-1}$,
\begin{equation}
 {\bm{p}}^{\ast} = {\mathop{\arg \min}_{\bm{p}}} \sum
 D_{\lambda}[{\bm{q}}_i:{\bm{p}}].
\end{equation}

The center is obtained from
\begin{equation}
 \partial_{\bm{p}} \sum_i D_{\lambda}
 \left[{\bm{q}}_i :{\bm{p}} \right] = 0,
\end{equation}
which yields the equation to give ${\bm{p}}^{\ast}$
\begin{equation}
 \sum {\rm{\bf G}} \left({\bm{q}}_i, {\bm{p}}^{\ast}\right)
 \left({\bm{q}}_i - {\bm{p}}^{\ast}\right)=0,
\end{equation}
where
\begin{equation}
 {\rm{\bf G}}({\bm{q}}, {\bm{p}}) = \nabla_{\bm{p}} \nabla_{\bm{p}}
 \varphi_{\lambda}({\bm{q}}, {\bm{p}})
\end{equation}

It is known that the mean (center) of two Gaussian distributions $N
\left(\mu_1, \sigma^2_1 \right)$ and $N \left(\mu_2, \sigma^2_2 \right)$
over the real line $X={\rm{\bf R}}$ is Gaussian $N \left( \frac{\mu_1+
\mu_2}2, \frac{\left(\sigma_1+ \sigma_2 \right)^2}4 \right)$, when we
use the square of the Wasserstein distance $W^2_2$ with the cost
function $|x_1-x_2|^2$. It would be interesting to see how the center changes depending on $\lambda$
based on $D_{\lambda}[{\bm{p}}:{\bm{q}}]$.

We consider the center of two Gaussian distributions ${\bm{q}}_1$ and ${\bm{q}}_2$, defined by
\begin{equation}
 {\bm{\eta}}_p ={\mathop{\arg\min}_{\bm{p}}} \sum D_{\lambda}
 \left[{\bm{p}}:{\bm{q}}_i \right].
\end{equation}
When $\lambda \rightarrow 0$ and $\lambda \rightarrow \infty$, we have
\begin{eqnarray}
 \lambda \rightarrow \infty &:& 
 \sigma^2_p = \frac{2 \sigma^2_{q_1}\sigma^2_{q_2}}{\sigma^2_{q_1}+
 \sigma^2_{q_2}},
  \quad
 \mu_p = \frac{\sigma^2_{q_2}\mu_{q_1}+
 \sigma^2_{q_1}\mu_{q_2}}{\sigma^2_{q_1}+ \sigma^2_{q2}}, \\
 \lambda \rightarrow 0 &:&
 \sigma_p = \frac{2 \sigma_{q_1} \sigma_{q_2}}{\sigma_{q_1}+
 \sigma_{q_2}}, \quad \mu_p = 
 \frac{\sigma_{q_2}\mu_{q_1} + \sigma_{q_1} \mu_{q_2}}{\sigma_{q_1}+
 \sigma_{q_2}}. 
\end{eqnarray}
However, if we use $C_{\lambda}$ instead of $D_{\lambda}$ the centers are
\begin{eqnarray}
 \lambda \rightarrow \infty &:& \sigma_p = \lambda, \\
 \lambda \rightarrow 0 &:& \sigma_p = 
  \frac{\sigma_{q_1}+ \sigma_{q_2}}2,
\end{eqnarray}
which are not reasonable for large $\lambda$.

\subsection{Statistical estimation} 
Let us consider a statistical model $M$,
\begin{equation}
 M = \left\{ p({\bm{x}}, {\bm{\xi}}) \right\}
\end{equation}
parameterized by $\bm{\xi}$. An interesting example is the set of distributions over $X= \left\{(0, 1)^n \right\}$, where ${\bm{x}}$ is a vector random variable defined on the $n$-cube $X$, where $\bm{x}$ is a vector random variable defined on the $n$-cube $X$. 

The Boltzmann machine $M$ is its submodel, consisting of probability distributions
which do not include higher-order interaction terms of random variables
$x_i$,
\begin{equation}
 p({\bm{x}}) = \exp \left\{
 \sum b_i x_i + \sum_{i<j} w_{ij} x_i x_j - \psi
 \right\}.
\end{equation} The transportation cost is
\begin{equation}
 {\rm{\bf m}}({\bm{x}}, {\bm{y}}) = \sum_i \left|x_i-y_i \right|,
\end{equation}
which is the Hamming distance \cite{Montavon2015}.

Let $\hat{\bm{q}}= \hat{\bm{q}}({\bm{x}})$ be an observed empirical
distribution. Then, $D_{\lambda}$-estimator ${\bm{p}}^{\ast}=
{\bm{p}}^{\ast}({\bm{x}}, {\bm{\xi}}^{\ast}) \in M$ is defined by
\begin{equation}
 {\bm{p}} \left({\bm{x}}, {\bm{\xi}}^{\ast}\right) =
 {\mathop{\arg \min}_{\bm{\xi}}} D_{\lambda}
  \left[\hat{\bm{q}} : p({\bm{x}}, {\bm{\xi}}) \right]. \label{eq:D_lam}
\end{equation}
Differentiating $D_{\lambda}$ with respect to  $\bm{\xi}$, we obtain the following theorem:
\theorem \upshape The $\lambda$-estimator ${\bm{\xi}}^{\ast}$ satisfies
\begin{equation}
 \label{eq:am8420161216}
 {\rm{\bf G}} \left(\hat{\bm{q}}, {\bm{p}}\right)
 \left({\bm{p}}-\hat{\bm{q}}\right)
 \frac{\partial p({\bm{x}}, {\bm{\xi}}^{\ast})}{\partial {\bm{\xi}}} = 0.
\end{equation}

\subsection{Pattern classifier}
Let ${\bm{p}}_1$ and ${\bm{p}}_2$ be two prototype patterns of
categories $C_1$ and $C_2$. A separating hyper-submanifold of the two
categories is defined by the set of ${\bm{q}}$ that satisfy
\begin{equation}
 D_{\lambda} \left[{\bm{p}}_1 : {\bm{q}}\right] =
 D_{\lambda} \left[{\bm{p}}_2 : {\bm{q}}\right]
\end{equation}
or
\begin{equation}
 D_{\lambda} \left[{\bm{q}}:{\bm{p}}_1 \right] =
 D_{\lambda} \left[{\bm{q}}:{\bm{p}}_2 \right].
\end{equation}

It would be interesting to study the geometrical properties of the $\lambda$-separating hyperplanes (Fig. 2).

\begin{figure}[h]
 \begin{center}
  \includegraphics[width=60mm]{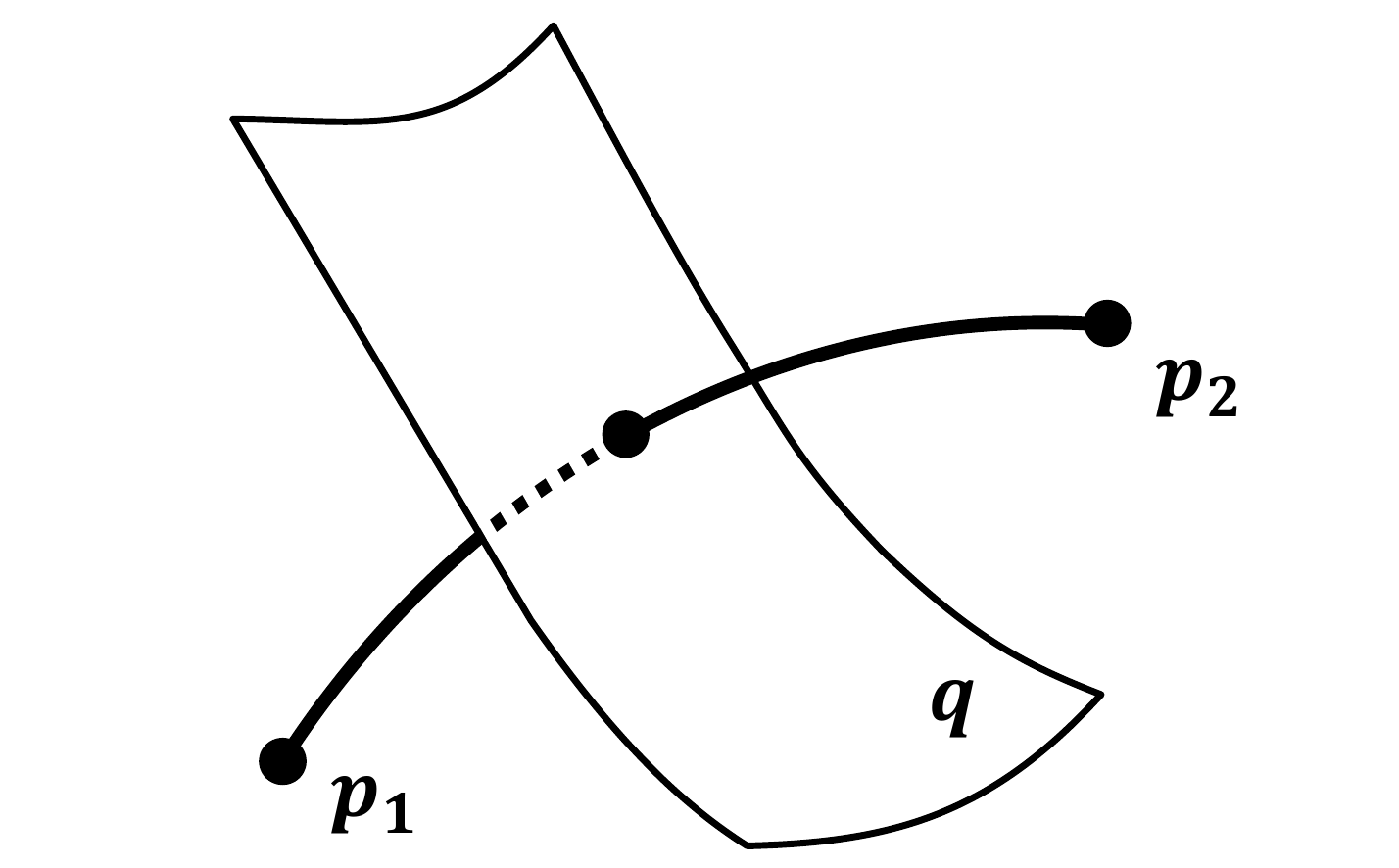}
 \caption{$\lambda$-separating hyperplane}
 \label{fig2}
\end{center}
\end{figure}

\section{Information Geometry of Transportation Plans}
We provide a general framework of the transportation plans from the viewpoint of information geometry. The manifold of all transportation plans is a probability simplex $M = S_{n^2-1}$ consisting of all the joint probability distributions ${\rm{\bf P}}$ over $X \times X$. It is dually flat, where $m$-coordinates are $\eta_{ij}= P_{ij}$, from which $P_{nn}$ is determined.
\begin{equation}
 \sum P_{ij} =1.
\end{equation}
The corresponding $e$-coordinates are $\log P_{ij}$, normalized by $P_{nn}$ as
\begin{equation}
 \theta^{ij}= \log \frac{P_{ij}}{P_{nn}}.
\end{equation}

We considered three problems in $M=S_{n^2-1}$, when the cost matrix ${\rm{\bf M}}= \left(m_{ij}\right)$ is given.

\ \ 

\noindent
{\textbf{1) Free problem}}

Minimize the entropy-relaxed transportation cost $\varphi_{\lambda}({\rm{\bf P}})$ without any constraints on ${\rm{\bf P}}$. The solution is
\begin{equation}
 {\rm{\bf P}}^{\ast}_{\rm{free}} = \exp \left(
  -\frac{m_{ij}}{\lambda} -\frac{1+\lambda}{\lambda} \psi \right) 
 = c {\rm{\bf K}},
\end{equation}
where $c$ is a normalization constant. This clarifies the meaning of the matrix {\rm{\bf K}} (Eq. (\ref{eq:K})), i.e., {\rm{\bf K}} is the optimal transportation plan for the free problem. 

\ \ 

\noindent
{\textbf{2) Rate-distortion problem}}

\begin{figure}[h]
 \begin{center}
  \includegraphics[width=55mm]{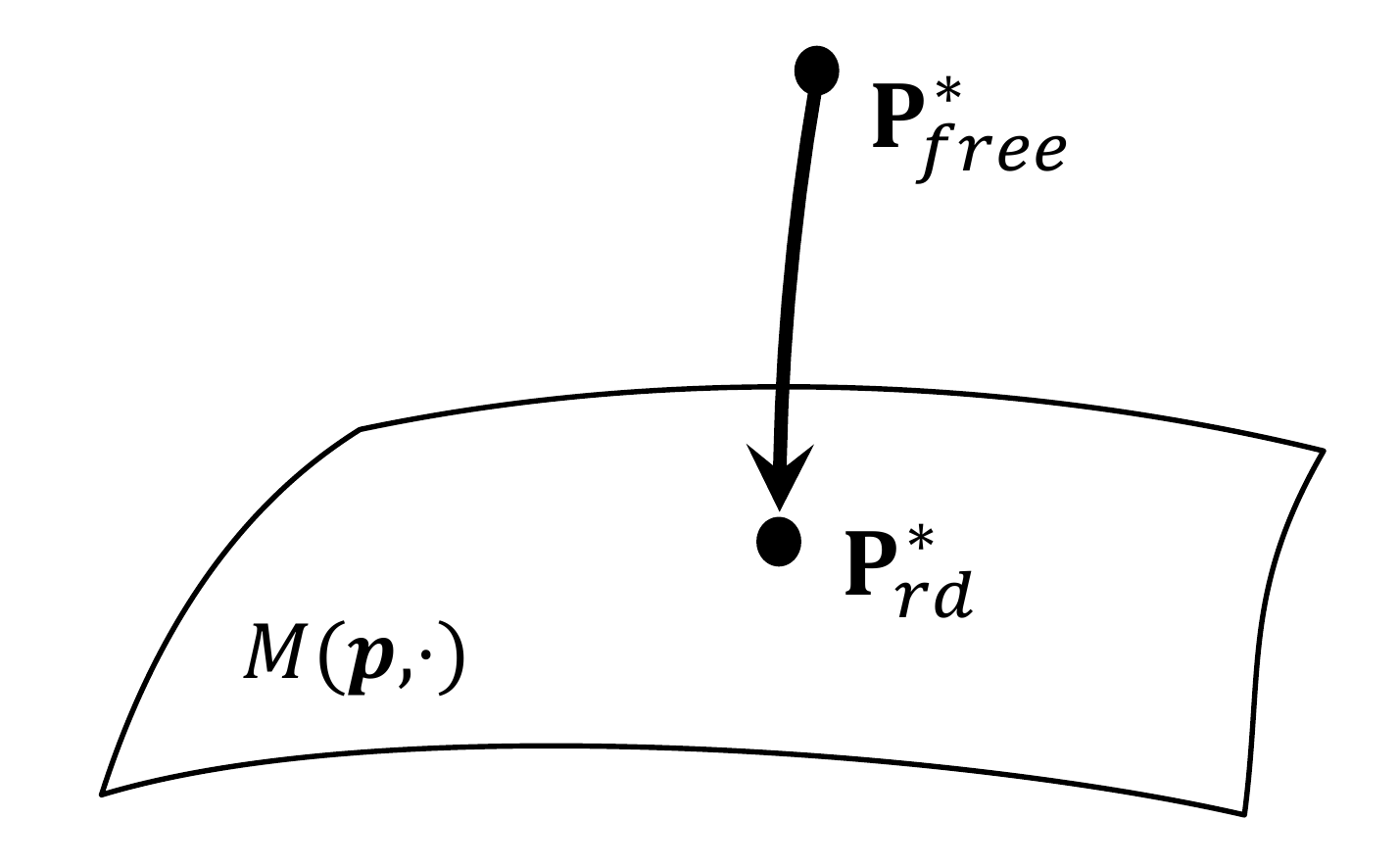}
 \caption{$e$-projection in the rate-distortion problem}
 \label{fig3}
 \end{center}
\end{figure}

\ \ 

We considered a communication channel in which ${\bm{p}}$ is a probability distribution on the sender’s terminals. The channel is noisy and $P_{ij}/p_i$ is the probability that $x_j$ is received when $x_i$ is sent. The costs $m_{ij}$ are regarded as the distortion of $x_i$ changing to $x_j$. The rate distortion-problem in information theory searches for ${\rm{\bf P}}$, which maximizes the mutual information of the sender and receiver under the constraint of distortion $\langle {\rm{\bf M}}, {\rm{\bf P}} \rangle$. The problem is formulated by maximizing $\varphi_{\lambda}(\rm{\bf P})$ under the sender’s constraint ${\bm{p}}$, where ${\bm{q}}$ is free (R. Belavkin, personal communication).

The optimal solution is given by
\begin{equation}
 {\rm{\bf P}}^{\ast}_{rd} = \left( c a_i K_{ij}\right),
\end{equation}
since ${\bm{q}}$ is free and ${\bm{\beta}}=0$ or $b_j=1$.  $a_i$ are determined from ${\bm{p}}$ such that the sender’s condition
\begin{equation}
 c \sum_j a_i K_{ij} = p_i
\end{equation}
is satisfied. Therefore, the dual parameters $a_i$ are given explicitly as
\begin{equation}
 \label{eq:am11720170126}
 c a_i = \frac{p_i}{\displaystyle{\sum_j K_{ij}}}.
\end{equation}

Let $M ({\bm{p}}, \cdot)$ be the set of plans that satisfy the sender’s condition 
\begin{equation}
 \label{eq:am11920170323}
 \sum_j P_{ij} = p_i.
\end{equation}
Then, we will see that ${\rm{\bf P}}^{\ast}_{rd}$ is the $e$-projection of ${\rm{\bf P}}^{\ast}_{\rm{free}}$ to $M({\bm{p}}, \cdot)$. The $e$-projection is explicitly given by Eq. (\ref{eq:am11720170126}) (Fig. 3).

\noindent
\\
{\textbf{3) Transportation problem}}
A transportation plan satisfies the sender’s and receiver’s conditions. Let $M(\cdot, {\bm{q}})$ be the set of plans that satisfies the receiver’s conditions \begin{equation}
 \sum_i P_{ij} = q_j.
\end{equation}
Then, the transportation problem searches for the plan that minimizes the entropy-relaxed cost in the subset
\begin{equation}
 M({\bm{p}}, {\bm{q}}) = M({\bm{p}}, \cdot) \cap M(\cdot, {\bm{q}}).
\end{equation}
Since the constraints Eqs. (\ref{eq:am11920170323}) and (102) are linear in the $m$-coordinates ${\rm{\bf P}}, M({\bm{p}}, \cdot)$, $M(\cdot, {\bm{q}})$ and $M({\bm{p}}, {\bm{q}})$ are $m$-flat submanifolds (Fig. 4).

\begin{figure}[h]
\begin{center}
  \includegraphics[width=70mm]{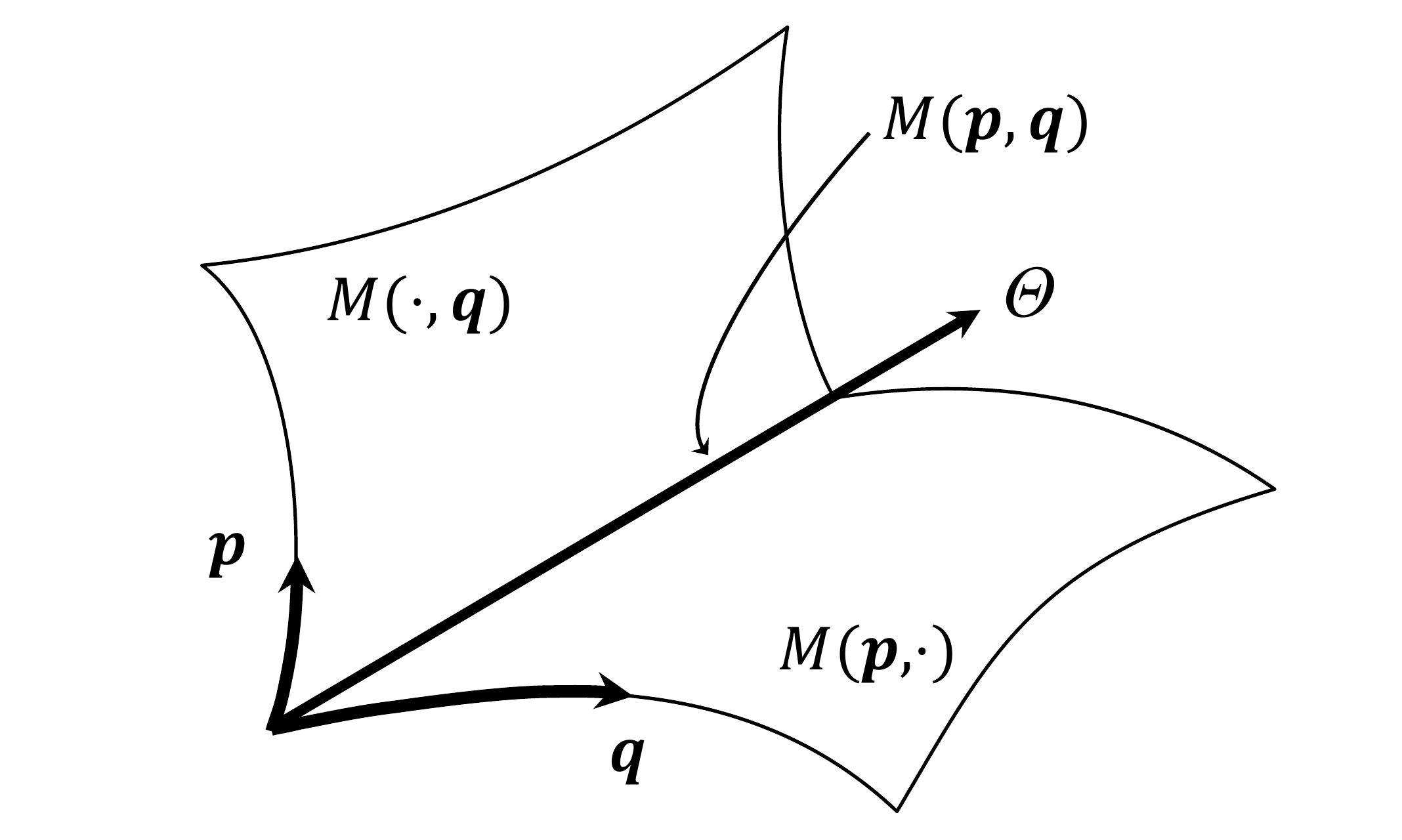}
 \caption{$m$-flat submanifolds in the transportation problem}
  \end{center}
 \label{fig4}
\end{figure}

Since $\bm{p}$ and $\bm{q}$ are fixed, $M({\bm{p}}, {\bm{q}})$ is of dimensions $(n-1)^2$, in which all the degrees of freedom represent mutual interactions between the sender and receiver. We define them by
\begin{equation}
 \label{eq:am7620161128}
  {\Theta}_{ij} = \log \frac{P_{ij}P_{nn}}{P_{in} P_{nj}},\quad
  i, j=1, \cdots, n-1.
\end{equation}
They vanish for ${\rm{\bf P}}_D={\bm{p}} \otimes {\bm{q}}$, as is easily seen Eq. (\ref{eq:am7620161128}).  Since ${\Theta}_{ij}$ are linear in $\log P_{ij}$, the submanifold $E \left({\Theta}_{ij}\right)$, in which ${\Theta}_{ij}$'s take fixed values but ${\bm{p}}$ and ${\bm{q}}$ are free, is an $2(n-1)$-dimensional $e$-flat submanifold.

We introduce mixed coordinates
\begin{equation}
 \Xi = \left({\bm{p}}, {\bm{q}}, {\Theta}_{ij} \right)
\end{equation}
such that the first $2(n-1)$ coordinates $({\bm{p}}, {\bm{q}})$ are the marginal distributions in the $m$-coordinates and the last $(n-1)^2$ coordinates $\Theta$ are interactions in the $e$-coordinates given in Eq. (\ref{eq:am7620161128}). Since the two complementary coordinates are orthogonal, we have orthogonal foliations of $S_{n^2-1}$ \cite{Amari2016} (Fig. 5).

\begin{figure}[h]
\begin{center}
  \includegraphics[width=65mm]{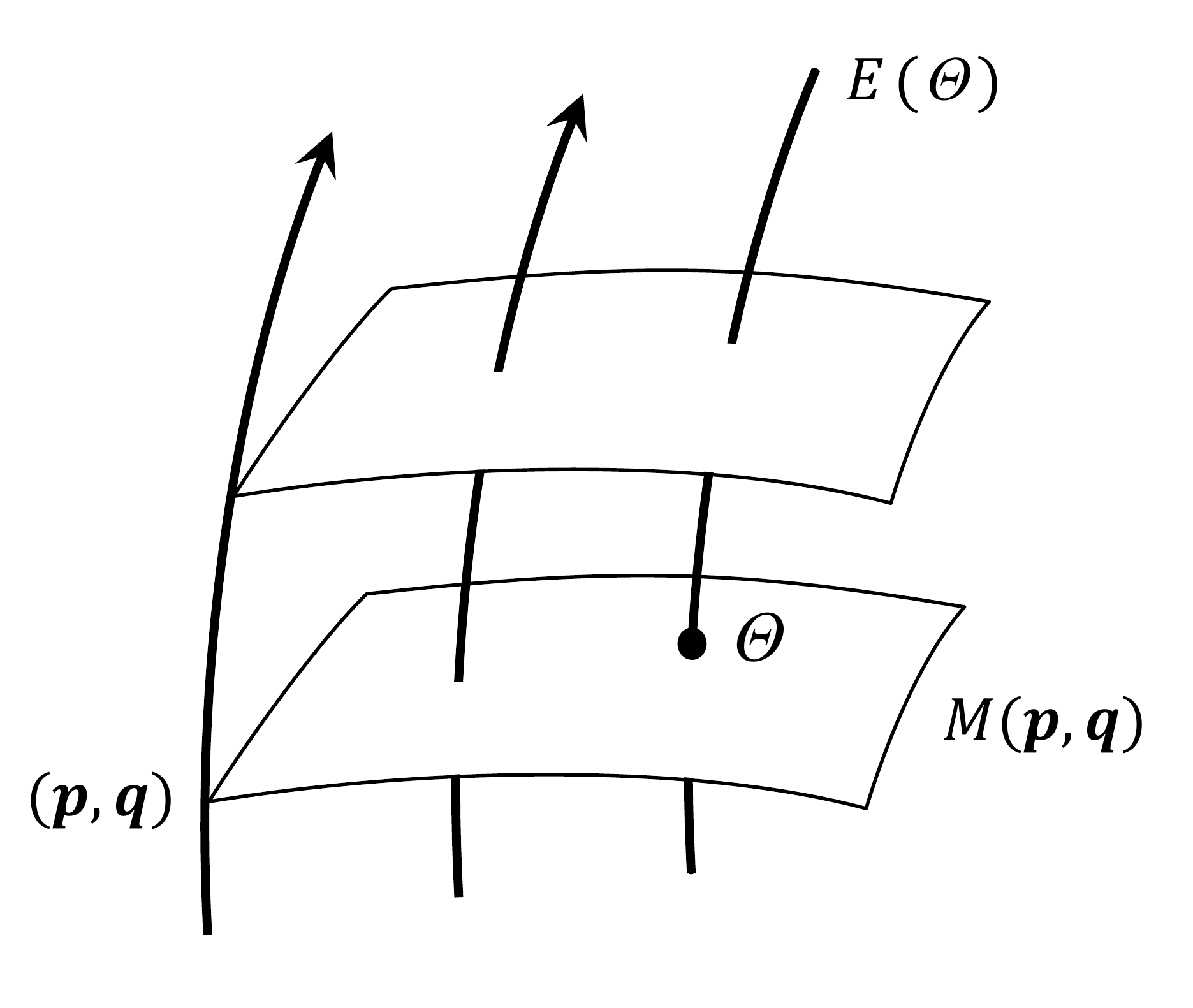}
 \caption{Orthogonal foliations of $S_{n^2-1}$ with the mixed coordinates}
  \end{center}
 \label{fig5}
\end{figure}

Given two vectors ${\bm{a}}= \left(a_i \right)$ and ${\bm{b}}= \left(b_j \right)$, we considered the following transformation of ${\rm{\bf P}}$,
\begin{equation}
 T_{\bm{a b}} {\rm{\bf P}} = \left(c a_i b_j {\rm{\bf P}}_{ij}\right),
\end{equation}
where $c$ is a constant determined from the normalization condition,
\begin{equation}
 c \sum_{i, j} a_i b_j P_{ij} = 1.
\end{equation}
$\Xi$ is the mixed coordinates of ${\rm{\bf P}}$ and $m$-flat submanifold $M({\bm{p}}, {\bm{q}})$, defined by fixing the first $2(n-1)$ coordinates, is orthogonal to $e$-flat submanifold $E\left(\Theta \right)$, defined by making the last $(n-1)^2$ coordinates equal to ${\Theta}_{ij}$. This is called the RAS transformation in the input-output analysis of economics.

\noindent 
\\
{\textbf{Lemma}}\quad For any ${\bm{a}}$, ${\bm{b}}$, transformation $T_{{\bm{a}}{\bm{b}}}$ does not change the interaction terms ${\Theta}_{ij}$. Moreover, the $e$-geodesic connecting ${\rm{\bf P}}$ and $T_{\bm{a b}}{\rm{\bf P}}$ is orthogonal to $M({\bm{p, q}})$.

\begin{proof} \upshape
By calculating the mixed coordinates of $T_{\bm{a b}}{\rm{\bf P}}$, we easily see that the ${\Theta}$-part does not change. Hence, the $e$-geodesic connecting ${\rm{\bf P}}$ and $T_{\bm{ab}}{\rm{\bf P}}$ is given, in terms of the mixed coordinates, by keeping the last part fixed while changing the first part. This is included in $E\left({\Theta}\right)$. Therefore, the geodesic is orthogonal to $M(\bm{p},\bm{q})$. 

\qed

Since the optimal solution is given by applying $T_{\bm{a b}}$ to ${\rm{\bf K}}$ such that the terminal conditions (Eq. (\ref{eq:am3})) are satisfied, we have the following theorem:  
\end{proof}

\theorem \upshape The optimal solution ${\rm{\bf P}}^{\ast}$ is given by $e$-projecting ${\rm{\bf K}}$ to $M({\bm{p}}, {\bm{q}})$.

\ \ 

\noindent {\textbf{4) Iterative Algorithm (Sinkhorn Algorithm) for obtaining $\bm{a}$ and $\bm{b}$}}
We need to calculate $\bm{a}$ and $\bm{b}$ when $\bm{p}$ and $\bm{q}$ are given for obtaining the optimal transportation plan. The Sinkhorn algorithm is well known for this purpose [11]. It is an iterative algorithm for obtaining the $e$-projection of ${\rm{\bf K}}$ to $M({\bm{p}}, {\bm{q}})$. 

Let $T_{A \cdot}$ be the $e$-projection of ${\rm{\bf P}}$ to $M({\bm{p}}, \cdot)$ and let $T_{\cdot B}$ be the $e$-projection to $M(\cdot, {\bm{q}})$. From the Pythagorean theorem, we have
\begin{equation}
 KL \left[T_{A \cdot}{\rm{\bf P}}:{\rm{\bf P}}\right]
 + KL \left[{\rm{\bf P}}^{\ast}: T_{A \cdot}{\rm{\bf P}}\right]
 = KL \left[{\rm{\bf P}}^{\ast}:{\rm{\bf P}}\right],
\end{equation}
where ${\rm{\bf P}}^{\ast}= T_{\bm{ab}}{\rm{\bf P}}$ is the optimal solution; that is, the $e$-projection of ${\rm{\bf K}}$ to $M({\bm{p}}, {\bm{q}})$. Hence, we have
\begin{equation}
 KL \left[{\rm{\bf P}}^{\ast}: T_{A \cdot} {\rm{\bf P}} \right]
 \le KL \left[{\rm{\bf P}}^{\ast}:{\rm{\bf P}} \right]
\end{equation}
and the equality holds when and only when ${\rm{\bf P}} \in M({\bm{p}}, \cdot)$. The $e$-projection of ${\rm{\bf P}}$ decreases the dual KL-divergence to ${\rm{\bf P}}^{\ast}$. The same property holds for the $e$-projection to $M(\cdot, {\bm{q}})$. The iterative $e$-projections of ${\rm{\bf K}}$ to $M({\bm{p}}, \cdot)$ and $M(\cdot, {\bm{q}})$ converges to the optimal solution ${\rm{\bf P}}^{\ast}$.

It is difficult to have an explicit expression of the $e$-projection of ${\rm{\bf P}}$ to $M({\bm{p}}, {\bm{q}})$, but those of $e$-projections to $M({\bm{p}}, \cdot)$ and $M(\cdot, {\bm{q}})$ are easily obtained. The $e$-projection of ${\rm{\bf P}}$ to $M({\bm{p}}, \cdot)$ is given by
\begin{equation}
 T_{A \cdot}{\rm{\bf P}} = \left(a_i P_{ij}\right),
\end{equation}
where ${\bm{a}}$ is given explicitly by
\begin{equation}
 a_i = \frac{p_i}{\sum_j P_{ij}}.
\end{equation}
Similarly, the e-projection to $M(\cdot, {\bm{q}})$ is given by
\begin{equation}
 T_{\cdot B}{\rm{\bf P}} = \left(b_j P_{ij}\right),
\end{equation}
with
\begin{equation}
 b_j = \frac{q_j}{\sum_i P_{ij}}.
\end{equation}
Therefore, the iterative algorithm, which is known as the Sinkhorn Algorithm \cite{Sinkhorn1964,Cuturi2013} of $e$-projection from
${\rm{\bf K}}$ is formulated as follows: 

\ \ 

\textbf{Iterative $e$-projection algorithm}
\begin{enumerate}
 \item[1.] Begin with ${\rm{\bf P}}_0={\rm{\bf K}}$.
 \item[2.] For $t=0, 1, 2, \cdots$, $e$-project $P_{2t}$ to $M({\bm{p}},
	   \cdot)$ to obtain
\begin{equation}
 {\rm{\bf P}}_{2t+1} = T_{A \cdot} {\rm{\bf P}}_{2t}.
\end{equation}
 \item[3.] To obtain ${\rm{\bf P}}_{2t+2}$, $e$-project ${\rm{\bf
	   P}}_{2t+1}$ to $M(\cdot, {\bm{q}})$,
\begin{equation}
 {\rm{\bf P}}_{2t+2} = T_{\cdot B} {\rm{\bf P}}_{2t+1}.
\end{equation}
 \item[4.] Repeat until convergence.
\end{enumerate}

\begin{figure}[h]
\begin{center}
  \includegraphics[width=70mm]{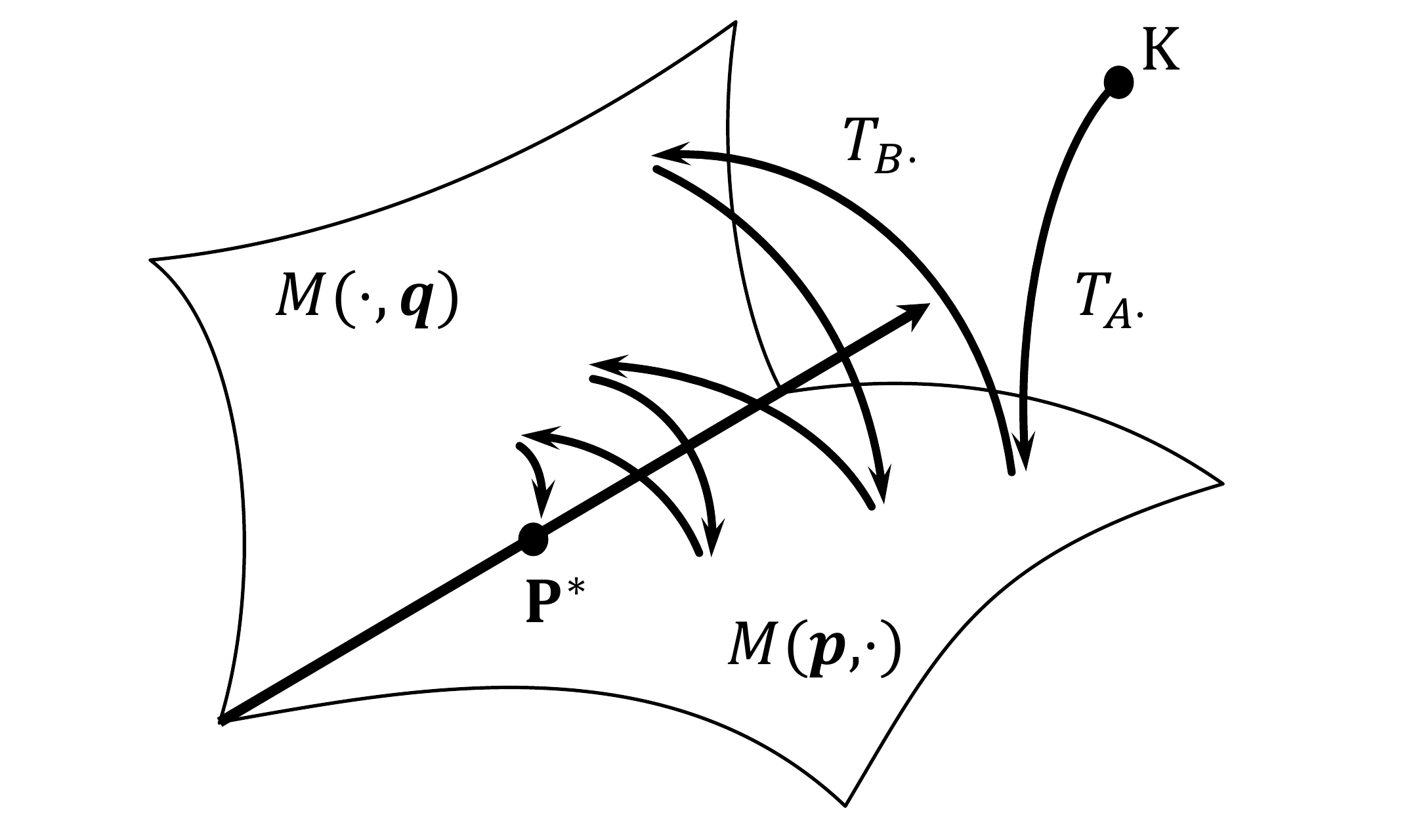}
 \caption{Sinkhorn algorithm as iterative $e$-projections}
  \end{center}
 \label{fig6}
\end{figure}

Fig. 6 schematically illustrates the iterative $e$-projection algorithm for finding the optimal solution ${\rm{\bf P}}^{\ast}$.

\section{Conclusions and Additional Remarks}

We elucidated the geometry of optimal transportation plans and introduced a one-parameter family of divergences in the probability simplex which connects the Wasserstein distance and KL-divergence. A one-parameter family of Riemannian metrics and dually coupled affine connections were introduced in $S_{n-1}$, although they are not dually flat in general. We uncovered a new way of studying the geometry of probability distributions. Future studies should examine the properties of the $\lambda$-divergence and apply these to various problems. We touch upon some related problems below. \\

\noindent {\textbf{1. Uniqueness of the optimal plan}} \\
The original Wasserstein distance is obtained by solving a linear programming problem.  Hence, the solution is not unique in some cases and is not necessarily a continuous function of ${\rm{\bf M}}$. However, the entropy-constrained solution is unique and continuous with respect to ${\rm{\bf M}}$ \cite{Cuturi2013}. While $\varphi_{\lambda}({\bm{p}}, {\bm{q}})$ converges to $\varphi_0({\bm{p}}, {\bm{q}})$ as $\lambda \rightarrow	0$, $\varphi_0({\bm{p}}, {\bm{q}})$ is not necessarily differentiable.

\noindent {\textbf{2. Integrated information theory of consciousness}} \\ 
Given a joint probability distribution ${\rm{\bf P}}$, the amount of integrated information is measured by the amount of interactions of information among different terminals. We used a disconnected model in which no information is transferred through branches connecting different terminals. The geometric measure of integrated information theory is given by the KL-divergence from ${\rm{\bf P}}$ to the submanifold of disconnected models \cite{Oizumi2016,Amari2017arXiv}. However, the Wasserstein divergence can be considered as such a measure when the cost of transferring information through different terminals depends on the physical positions of the terminals \cite{Oizumi2014}. We can use the entropy-constrained divergence $D_{\lambda}$ to define the amount of information integration.

\noindent {\textbf{3. $f$-divergence}} \\
We used the KL-divergence in a dually flat manifold for defining $D_{\lambda}$. It is possible to use any other divergences, for example, the $f$-divergence instead of KL-divergence. We would obtain similar results. 

\noindent {\textbf{4. $q$-entropy}} \\
Muzellec et al. used the $\alpha$-entropy (Tsallis $q$-entropy) instead of the Shannon entropy for regularization \cite{Muzellec2016}. This yields the $q$-entropy-relaxed framework.

\noindent {\textbf{5. Comparison of $C_{\lambda}$ and $D_{\lambda}$}} \\
Although $D_{\lambda}$ satisfies the criterion of a divergence, it might differ considerably from the original $C_{\lambda}$. In particular, when $C_\lambda({\bm{p}}, {\bm{q}})$ includes a piecewise linear term such as $\sum d_i|p_i - q_i|$ for constant $d_i$, $D_\lambda$ defined in Eq. (\ref{eq:Bregman_C}) eliminates this term. When this term is important, we can use $\{ C_\lambda({\bm{p}}, {\bm{q}}) \}^2$ instead of $C_\lambda({\bm{p}}, {\bm{q}})$ for defining a new divergence $D_\lambda$ in Eq. (\ref{eq:Bregman_C}). In our accompanying paper \cite{Amari2017}, we define a new type of divergence that retains the properties of $C_{\lambda}$ and is closer to $C_{\lambda}$.



\section*{Appendix: The Proof of Example 2}
\renewcommand{\theequation}{A.\arabic{equation}}
\setcounter{equation}{0}

Let us assume that functions $a(x)$ and $b(y)$ are constrained into Gaussian distributions:     
$a(x) = N
\left(\tilde{\mu}, \tilde{\sigma}^2 \right)$, $b(y) = N
\left(\tilde{\mu}', \tilde{\sigma}'^2 \right)$.
This means that the optimal plan $ P^{\ast}(x, y)$ is also given by a Gaussian distribution $N ({\boldsymbol \mu},\Sigma)$.
The marginal distributions $p$ and $q$ require the mean value of the optimal plan to become   
\begin{equation}
{\boldsymbol \mu} = [\mu_p \ \ \mu_q]^T.
\end{equation}
It is also necessary for the diagonal part of the covariance matrix to become 
\begin{eqnarray}
&& \Sigma_{11} = \sigma_p^2, \\
&& \Sigma_{22} = \sigma_q^2.
\end{eqnarray}
Because the entropy-relaxed optimal transport is given by Eq. (\ref{GaussianP}), $\Sigma$ is composed of $\tilde{\sigma}^2$ and $\tilde{\sigma}'^2$
as follows: 
\begin{eqnarray}
&& \Sigma_{11} = \frac{\tilde{\sigma}^2 (2 \tilde{\sigma}'^2 +\lambda)}{ 2(\tilde{\sigma}^2+\tilde{\sigma}'^2) + \lambda}, \\
&& \Sigma_{22} =  \frac{\tilde{\sigma}'^2 (2 \tilde{\sigma}^2 +\lambda)}{ 2(\tilde{\sigma}^2+\tilde{\sigma}'^2) + \lambda}. 
\end{eqnarray}
Solving Eqs. (A.4,5) under the conditions given in Eqs. (A.2,3), we have
\begin{eqnarray}
&& \tilde{\sigma}^2  = \left \{ \frac{1}{2\sigma_p^2} (1+\sqrt{1+X}) -\frac{2}{\lambda} \right \}^{-1}, \\
&& \tilde{\sigma}'^2 =  \left \{ \frac{1}{2\sigma_q^2} (1+\sqrt{1+X}) -\frac{2}{\lambda} \right \}^{-1},  \\
&& \mbox{where} \quad X= \frac{16 \sigma^2_p \sigma^2_q}{\lambda^2}.
\end{eqnarray}
Substituting the mean (Eq. (A.1)) and variances (Eqs. (A.6,7)) into the definition of the cost (Eq. (\ref{eq:Cuturi_cost})), after straightforward calculations, we get Eq. (\ref{GaussianC}).
In general, the $\eta$ coordinates of the Gaussian distribution $q$ are given by $(\eta_1,  \eta_2) = (\mu_q, \mu_q^2+ \sigma_q^2)$. 
After differentiating $C_{\lambda}(p,q)$ with the $\eta$ coordinates and substituting them into Eq. (\ref{eq:Bregman_C}), we get Eq. (\ref{GaussianD}).

\end{document}